\documentclass[11pt]{article}
\usepackage{epsfig}
\usepackage{psfrag}
\usepackage{amsmath,amssymb,amsthm,enumerate}

\usepackage{cite}
\begin{document}
%

%
\newtheorem{thm}{Theorem}
\newtheorem{defn}[thm]{Definition}
\newtheorem{cor}[thm]{Corollary}
\newtheorem{prop}[thm]{Proposition}
\newtheorem{lemma}[thm]{Lemma}
\newcommand{\beqn}{\begin{equation}}
\newcommand{\eeqn}{\end{equation}}
\newcommand{\nn}{\nonumber}
\newcommand{\cleq}{\preccurlyeq}
\newcommand{\cgeq}{\succcurlyeq}
\newcommand{\la}{{\langle}}
\newcommand{\ra}{{\rangle}}
\def\pa{{\partial}}
\def\R{{\mathbb R}}
\def\C{{\mathbb C}}
\def\ep{{\epsilon}}

\def\N{{\cal N}}
\def\vb{{\bar{v}}}
\def\phib{{ \bar{\phi}}}
\def\vbb{{\bar{\vb}}}
\def\ub{{\bar{u}}}
\def\wb{{\bar{w}}}
\def\wbb{{\bar{\wb}}}
\def\wt{{\tilde{w}}}
\def\wtt{{\tilde{\wt}}}
\def\Ft{{\tilde{F}}}

\numberwithin{thm}{section}
\numberwithin{equation}{section}

\title{Uniqueness for a hyperbolic inverse problem with angular control on the coefficients}

\author{
Rakesh\thanks{Partially supported by an NSF Grant DMS-0907909} \\
Department of Mathematical Sciences\\
University of Delaware\\
Newark, DE 19716, USA\\
~\\
Email: rakesh@math.udel.edu\\
\and
Paul Sacks\\
Department of Mathematics\\
Iowa State University\\
Ames, IA 50011\\
~\\
Email: psacks@iastate.edu
 }
\date{December 15, 2010}
\maketitle



\begin{abstract}
Suppose $q_i(x)$, $i=1,2$ are smooth functions on $\R^3$ and $U_i(x,t)$ the solutions of the initial value problem
\begin{gather*}
\pa_t^2 U_i- \Delta U_i - q_i(x) U_i = \delta(x,t), \qquad (x,t) \in \R^3 \times \R
\\
U_i(x,t) =0, \qquad \text{for} ~ t<0.
\end{gather*}
Pick $R,T$ so that $0 < R < T$ and let $C$ be the vertical cylinder $\{ (x,t) \, : |x|=R, ~ R \leq t \leq T \}$.
We show that if $(U_1, U_{1r}) = (U_2, U_{2r})$ on $C$ then $q_1 = q_2$ on the annular region $R \leq |x| \leq (R+T)/2$
provided there is a $\gamma>0$, independent of $r$, so that
\[
 \int_{|x|=r} | \Delta_S (q_1 - q_2)|^2 \, dS_x \leq \gamma
 \int_{|x|=r}  |q_1 - q_2|^2 \, dS_x,
 \qquad \forall r \in [R, (R+T)/2].
\]
Here $\Delta_S$ is the spherical Laplacian on $|x|=r$. 
\end{abstract}


\section{Introduction}

Our goal is the study of a formally determined inverse problem for a hyperbolic PDE. Consider an acoustic medium, occupying the region $\R^3$, excited by an impulsive point source and the response of the medium is measured for a certain time period at receivers placed on a sphere surrounding the source. We study the question of recovering the acoustic property of the medium from this measurement.

Let $q(x)$ be a smooth function on $\R^3$ and $U(x,t)$ the solution of the initial value problem
\begin{gather}
U_{tt} - \Delta U - q(x)U = 8 \pi \delta(x,t), \qquad (x,t) \in
\R^3 \times \R,
\label{eq:Ude}
\\
U = 0, \qquad t< 0.
\label{eq:Uic}
\end{gather}
Using the progressing wave expansion one may show that
\beqn
U(x,t) = 2 \frac{ \delta (t-|x|)}{|x|}  +  u(x,t) H(t-|x|),
\label{eq:Uu}
\eeqn
where $u(x,t)$ is the solution of the Goursat problem
\begin{gather}
u_{tt} - \Delta u - q(x) u = 0, \qquad (x,t) \in \R^3, ~ t \geq |x|,
\label{eq:ude}
\\
u(x,|x|) = \int_0^1 q( \sigma x) \, d \sigma.
\label{eq:uic}
\end{gather}
The well posedness of the above Goursat problem is proved in \cite{Rom87} and improved in \cite{Rom09}, though the result is not optimal; \cite{Rom87} has suggestions for obtaining better results and we will address them elsewhere. For completeness we restate the well posedness result.
\begin{thm}[{See \cite{Rom87} and \cite{Rom09}}]\label{thm:existence}
Suppose $\rho>0$, and $q$ is a $C^{8}$ function on the ball $|x| \leq \rho$; then (\ref{eq:ude}), (\ref{eq:uic}) has a unique $C^2$ solution on the double conical region
$\{ (x,t) \in \R^3 \times \R \, : \, |x| \leq \rho, ~ |x| \leq t \leq 2\rho - |x| \}$. Further,
the $C^2$ norm of $u$, on this double conical region, approaches zero if the $C^8$ norm of $q$, on $|x| \leq \rho$, approaches zero. Also, if $q$ is smooth then so is $u$.
\end{thm}

\begin{figure}
 \begin{center}
  \epsfig{file=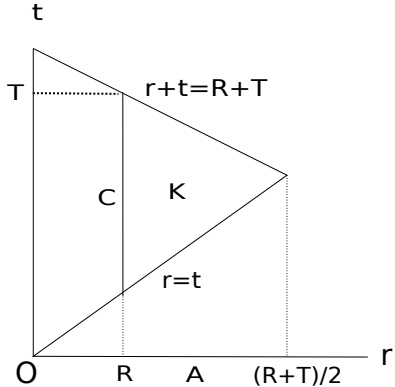, height=2in}
\end{center}
  \caption{Geometry of the problem}
\label{fig:vbasic}
\end{figure}

Below $P \cleq Q$ will mean that $P \leq C Q$ for some constant $C$.
Let $S$ denote the unit sphere centered at the origin. For any $0<R<T$, we define (see Figure \ref{fig:vbasic}) the annular region
\[
A := \{ x \in \R^3 \, : \, R \leq |x| \leq (R+T)/2 \},
\]
the space-time cylinder
\[
C = \{ (x,t) \in \R^3 \times \R \, : \, |x|=R, ~ R \leq t \leq T \},
\]
and
\[
K := \{ (x,t) \in \R^3 \times \R \, : \, R \leq |x| \leq (R+T)/2, ~ |x| \leq t \leq R+T-|x| \},
\]
a region bounded by $C$ and two light cones. 

In our model the source is at the origin, the receivers are on the sphere $|x|=R$ and the signals are measured up to time $T$.
Hence we define the forward map  
\[
F \, : \, q \mapsto (u|_C, u_r|_C)
\]
and our goal is to study the injectivity and the inversion of $F$. From the domain of dependence property of solutions of hyperbolic PDEs, it is clear that $F(q)$ is unaffected by changes in $q$
in the region $|x| \geq (R+T)/2$. Hence the best we can hope to do is recover $q$ on the ball $|x| \leq (R+T)/2$. 

If $q$ is spherically symmetric then the problem reduces to an inverse problem for the one dimensional wave equation. In this case, recovering $q$ on the region $R \leq |x| \leq (R+T)/2$, from $F(q)$, is done by the downward continuation method or the layer stripping method - see \cite{Sym83} and other references there. However, even in the spherically symmetric case (i.e.~the one dimensional case), recovering $q$ on $|x| \leq R$, from $F(q)$ is more difficult since the downward continuation scheme is not directly applicable. It is believed that uniqueness does not hold for this inverse problem if $T < 3R$ though explicit examples have not been constructed. If $T \geq 3R$, the question of recovering $q$ on $|x| \leq R$ from $(u,u_r)|_C$ was resolved by connecting this problem to one where the downward continuation method is applicable - see \cite{Rak00} and the references there. So it seems that in the general $q$ case, recovering $q$ over the region $|x| \leq R$ will be harder than recovering $q$ over the region $R \leq |x| \leq (R+T)/2$. 

\begin{figure}[h]
\begin{center}
\epsfig{file=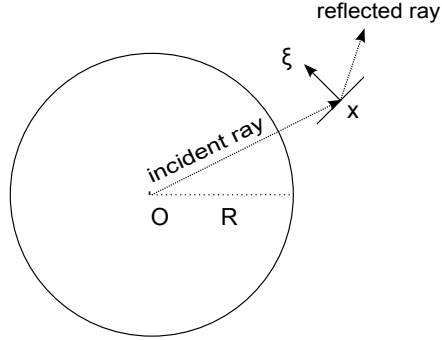, height=1.75in}
\end{center}
\caption{Reflection by a singularity in $q$}
\label{fig:wavefront}
\end{figure}

Our main result concerns the problem of recovering $q$ on $R \leq |x| \leq (R+T)/2$ from $(u,u_r)|_C$. The downward continuation method does not apply directly in higher space dimensions since the time-like Cauchy problem for hyperbolic PDEs is ill-posed in higher space dimensions. Further, an analysis of the linearized problem shows that there could be singularities in $q$ in certain directions, that is points in the wave front set of $q$, so that a signal emanating from the origin is reflected by this singularity in $q$, and the reflected signal never reaches the sphere $|x|=R$ where the receivers are located - see Figure \ref{fig:wavefront}. Hence there should not be any stability for this inverse problem, unless we restrict $q$ to a class of functions where singularities in $q$ of the above type are controlled. In \cite{SaSy85},  Sacks and Symes adapted the downward continuation method to apply to a slightly different inverse problem, with an impulsive planar source $\delta(z-t)$, with data measured on the hypersurface $z=0$, where $x=(y,z)$ with $y \in \R^2$ and $z \in \R$. They proved uniqueness for the linearized inverse problem when the unknown coefficient was restricted to the class of functions whose derivatives in the $y$ direction  were controlled by derivatives in the $z$ direction. Later Romanov showed the inversion methods for one dimensional problems could be used for the existence and reconstruction for the nonlinear version of the Sacks and Symes inverse problem provided $q(y,z)$ lies in the class of functions which are analytic in $y$ in a certain sense, that is strong restrictions are placed on the changes in $q$ in the $y$ direction - see \cite{Rom02} for details. We apply the technique in \cite{SaSy85} to the uniqueness question for the problem of 
recovering $q$ on on $R \leq |x| \leq (R+T)/2$ from $(u,u_r)|_C$; we will have to impose restrictions on the angular derivatives of $q$.

For any $x \in \R^3$ we define $r=|x|$ and for $x \neq 0$ we define $\theta = x/r \in S$; hence $x=r \theta$.
Define the radial vector field $\pa_r = r^{-1} x \cdot \nabla$ and, for $1 \leq i < j \leq 3$, the angular vector fields
$\Omega_{ij} = x_i \pa_j - x_j \pa_i$.
\begin{defn}
Given $\gamma>0$, we define $Q_\gamma(R,T)$ to be the set of all $C^2$ functions $q(x)$ on the ball $|x| \leq (R+T)/2$ with
\[
\| p\|_{H^2(S_r)} + \| \pa_r p \|_{H^1(S_r)} 
\leq \gamma  \left ( \|p\|_{H^1(S_r)} + \|\pa_r p\|_{L^2(S_r)} \right )
 \qquad  \forall r \in [R,(R+T)/2]
\]
where $p(x) = \int_0^{|x|} q( \sigma x/|x|) \, d \sigma$ and $S_r$ is the sphere $|x|=r$.
\end{defn}
So if $q$ is a smooth function on $|x| \leq (R+T)/2$ with
$\|p\|_{H^1(S_r)} + \|\pa_r p\|_{L^2(S_r)}  $ nonzero for every 
$r \in [R, (R+T)/2]$ then $q \in Q_\gamma$ where
\[
\gamma = \max_{ r \in [R, (R+T)/2] }
\dfrac{ \| p\|_{H^2(S_r)} + \| \pa_r p \|_{H^1(S_r)}  }
{\|p\|_{H^1(S_r)} + \|\pa_r p\|_{L^2(S_r)} }.
\]
Noting that $\pa_r p=q$, using Garding's inequality on a 
sphere\footnote{The Euclidean version is (6.8) on page 66 of \cite{La85}. Using a partition of unity argument and the Euclidean version, one may show that $ \|q\|_{H^2(S_r)} \leq C_r \| \Delta_S q\|_{L^2(S_r)}$ with $C_r$ bounded if $r$ is in a closed interval not containing $0$.}
, one may show that $q \in Q_{\gamma^*}$ for some $\gamma^* > 0$ if there is a $\gamma>0$ so that 
\[
\| \Delta_S q \|_{L^2(S_r)} \leq \gamma \|q\|_{L^2(S_r)},
\qquad \forall r \in [R, (R+T)/2]
\]
where $\Delta_S$ is the Laplacian on $S_r$.
In particular, if $q$ is a finite linear combination of the spherical harmonics with coefficients dependent on $r$ then $q \in Q_\gamma(R,T)$ for some $\gamma>0$.

In section \ref{sec:reflection} we prove the following injectivity result using the ideas in \cite{SaSy85}. 
\begin{thm}\label{thm:reflection}
Suppose $0<R<T$ and $q_1, q_2$ are $C^8$ functions on $\R^3$.
If $F(q_1) = F(q_2)$ and $q_1-q_2 \in Q_\gamma(R,T)$ for some $\gamma>0$ then $q_1=q_2$ on $R \leq |x| \leq (R+T)/2$.
\end{thm}

One may tackle the problem dealt with in Theorem \ref{thm:reflection}  using Carleman estimates also and one obtains a result which is stronger in some aspects and weaker in others. Using Carleman estimates one can prove uniqueness under slightly less stringent conditions on $q$ - one needs controls on the $L^2$ norms of only the first order angular derivatives of $p$ in terms of the $L^2$ norm of $p$, instead of on the second order angular derivatives required in Theorem \ref{thm:reflection}. However, the price one pays is that the $\gamma$ cannot be arbitrary but is determined by $R,T$; further $R$ cannot be arbitrary, but must satisfy $R>T/2$ and uniqueness is proved only for the values of $q$ in an annular region $R \leq |x| \leq R^*$ for some $R^* < (R+T)/2$. This work will appear elsewhere.

From Theorem \ref{thm:reflection} we can easily derive the following interesting corollary.
\begin{cor}\label{cor:injectivity}
Suppose $0<T$ and $q_1, q_2$ are smooth functions on $\R^3$ which vanish in a neighborhood of the origin. If $u_1$ and $u_2$ agree to infinite order on the line $\{(x=0,t) \, : \, 0 \leq t \leq T \}$ and $q_1-q_2 \in Q_\gamma(0,T)$ for some $\gamma>0$, then 
$q_1 = q_2$ on $|x| \leq T$.
\end{cor}
We give a short proof of the corollary. If $q_1=q_2 =0$ in some small neighborhood of the origin then the difference $u=u_1-u_2$ satisfies the standard homogeneous wave equation in a semi-cylindrical region 
\beqn
\{ (x,t) \in \R^3 \times R \, : \, |x| \leq \delta, ~ |x| \leq t \leq T-|x| \},
\label{eq:semicyl}
\eeqn
for some $\delta>0$.
Now, from the hypothesis, we have $u$ is zero to infinite order on the segment of the $t$ axis consisting of $0 \leq t \leq T$.
Then by Lebeau's unique continuation result in \cite{Leb99} we have $u=0$ in the semi-cylindrical region given in (\ref{eq:semicyl}).
Hence $u$ and $u_r$ are zero on the cylinder 
\[
\{ (x,t) \in \R^3 \times R \, : \, |x| = \delta, ~ \delta \leq t \leq T-\delta \}.
\]
The corollary follows from Theorem (\ref{thm:reflection})  if the $R$ and $T$ in Theorem \ref{thm:reflection} are taken to be $\delta$ and $T - \delta$ respectively.

We also have a uniqueness result for the linearized version of the inverse problem considered in Theorem \ref{thm:reflection}; the result is for a linearization about a radial background.
\begin{thm}\label{thm:linearized}
Suppose $q_b(r)$ is a function on $[0,\infty)$ so that $q_b(|x|)$ is a smooth function on $\R^3$; further suppose $u_b(r,t)$ is the solution of (\ref{eq:ude}), (\ref{eq:uic}) when $q(x)$ is replaced by $q_b(|x|)$. Let $q(x)$ be a smooth function on $\R^3$ and $u(x,t)$ the solution of the Goursat problem
\begin{gather}
u_{tt} - \Delta u - q_b u = q u_b, \qquad t \geq |x|,
\label{eq:ulinde}
\\
u(x,|x|) = \int_0^1 q( \sigma x) \, d \sigma.
\label{eq:ulinic}
\end{gather}
If $(u,u_r)|_C=0$ then $q=0$ on the region $R \leq |x| \leq (R+T)/2$.
\end{thm}
This theorem holds with less regular $q_b$ and $q$; what is needed is enough regularity so that the spherical harmonic expansions of $q$, $q_b$ and $u_b$ converge in the $C^2$ norm.

\begin{figure}[h]
\begin{center}
\epsfig{file=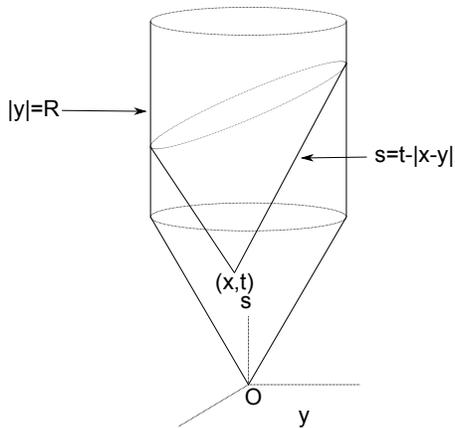, height=2.2in}
\end{center}
\caption{Kirchhoff's Formula}
\label{fig:kirk}
\end{figure}

We next focus on the  problem of recovering $q$ on the region $|x| \leq R$ from $(u, u_r)|_C$ when $T\geq 3R$. 
The linearized problem about the $q=0$ background, consisting of recovering $q$ from $(u, u_r)|_C$, where $u(x,t)$ is the solution of
the Goursat problem
\begin{gather*}
u_{tt} - \Box u = 0, \qquad t \geq |x|,
\\
u(x,|x|) = \int_0^1 q( \sigma x) \, d\sigma.
\end{gather*}
As observed by Romanov, since $T \geq 3R$, we may recover $q$ from  $(u, u_r)|_C$ fairly quickly. In fact, from Kirchhoff's formula (see \cite{Fr75}) expressing the solution of the wave equation in terms of the Cauchy data on $C$, we have
\[
u(x,t) = \int_{|y-x|=R} \frac{ u_r(y, t+|x-y|)}{|x-y|} +
\left (  \frac{ u(y, t+|x-y|)}{|x-y|^2}  +  \frac{ u_t(y, t+|x-y|)}{|x-y|}  \right ) \, \frac{ (y-x) \cdot y}{|x-y|} \, dS_y.
\]
for all $(x,t)$ with $|x| \leq t \leq R$ - see Figure \ref{fig:kirk}. In particular we can express $u(x,|x|)$ in terms of $(u, u_t, u_r)|_C$ and hence we can recover $q$.

For the {\em original nonlinear inverse problem} we show a partial uniqueness and stability result when one of the $q$ is small.
\begin{thm}\label{thm:transmission}
Suppose $0<3R<T$, $M>0$ and $q_i$, $i=1,2$ are $C^{8}$ functions on $|x| \leq (R+T)/2$ with $\|q_i\|_\infty \leq M$.
Let $u_i$ be the unique solution of (\ref{eq:ude}), (\ref{eq:uic}) with $q$ replaced by $q_i$; then there is a constant $\delta>0$ depending only on $R,T$ and $M$ so that if $\|q_2\|_\infty \leq \delta$ then
\beqn
\int_{|x| \leq R} |q_1-q_2|^2 \, dx \cleq 
\int_C |u_1-u_2|^2 + |\nabla (u_1 - u_2)|^2 + |(u_1 - u_2)_t|^2 \, dS_{x,t};
\label{eq:smallresult}
\eeqn
the constant in (\ref{eq:smallresult}) depending only on $R,T,M$.
\end{thm}
A weaker form of this result, requiring that $\|q_1\| \leq \delta$ also, was given in \cite{Rak93}; a result similar to this weaker version was also derived in \cite{RomYam99}. Later it was observed in \cite{LiYam07}, for a similar type of problem, that the above proofs go
through without the extra assumption that $\|q_1\| \leq \delta$. We give this short proof of Theorem \ref{thm:transmission}, in section \ref{sec:transmission}. However, the original nonlinear inverse problem remains unsolved.


\section{Proof of Theorem \ref{thm:reflection}}\label{sec:reflection}

\subsection{Preliminary observations}
We need the following observations in the proof.
For the angular vector fields we have $[\Omega_{ij}, \pa_r] = 0$, and $[\Omega_{ij}, \Omega_{kl}] =0$ if $\{i,j\}=\{k,l\}$ but $[\Omega_{ij}, \Omega_{ik}] = \Omega_{kj}$.  Also $|\nabla f|^2 = f_r^2 + r^{-2} \sum_{i<j} (\Omega_{ij} f)^2$ and if we define $\Omega =  \sum_{i<j} \Omega_{ij}^2$
then $\Delta = \pa_r^2 + 2r^{-1} \pa_r + r^{-2}\Omega$ and $[\Omega_{ij}, \Delta]=0$. 
Also,
for any $i \neq j$, since $ \Omega_{ij} f =  x_i \pa_j f - x_j \pa_i f = \pa_j( x_i f) - \pa_i( x_j f)$ and $x_j x_i - x_i x_j =0$,  by the divergence theorem, for any $0 < R_1 < R_2$ we have
\beqn
\int_{R_1 \leq |x| \leq R_2}  \Omega_{ij} f \, dx = 0.
\label{eq:divergence}
\eeqn
Applying (\ref{eq:divergence}) to the zeroth order homogeneous extension of a function $f$ on $S$, we conclude that for $C^1$ functions $f,g$ on $S$
\beqn
\int_S \Omega_{ij} f \, dS =0,
\qquad  
\int_S f \, \Omega_{ij} g \, dS = - \int_S g \, \Omega_{ij} f \, dS.
\label{eq:stokes}
\eeqn

For $i=1,2$ let $u_i$ be the solution of (\ref{eq:ude}), (\ref{eq:uic}) when $q=q_i$. Define $v_i(x,t)=ru_i(x,t)$,  $p_i(x) = r \int_0^1 q_i(\sigma x) \, d \sigma = \int_0^r q_i(\sigma \theta) \, d \sigma$. Define $v = v_1-v_2$, $q=q_1-q_2$ and
$p=p_1-p_2$. Then we have
\begin{gather}
v_{tt} - v_{rr} - \frac{1}{r^2} \Omega v - q_1v = q v_2, \qquad t \geq |x|
\label{eq:vde}
\\
v(x,|x|) = p(x).
\label{eq:vic}
\end{gather}
We are given that $(v, v_r)$ are zero on $C$ and we have to show that $q=0$ on $R \leq |x| \leq (R+T)/2$. Note that since $v=0$ on $C$, we have $p(x) = v(x,|x|) =0$ on $|x|=R$ and hence for $|x| \geq R$ we have $p(x) = \int_R^r q(\sigma \theta) \, d \sigma$ and hence to prove the theorem it will be enough to show that $p(x)=0$ on $R \leq |x| \leq (R+T)/2$.

We will attempt to carry out the proof which works in the one dimensional case. The limitations of this method when applied to the three dimensional case force the restrictions on $q$ in the statement of Theorem \ref{thm:reflection}. In the one dimensional case the angular terms are missing from (\ref{eq:vde}) so the roles of $r,t$ are reversible and one has sideways energy estimates which allow us to estimate the $H^1$ norm of $v$ on $t=|x|$ in terms of the norm of $v, v_r$ on $r=R$ and the $L^2$ norm of the RHS of (\ref{eq:vde}). The $H^1$ norm of $v$ on $t=|x|$ dominates the $L^2$ norm of $q$ on $A$ and the  $L^2$ norm of the RHS of the (\ref{eq:vde}) is dominated by $T-R$ times the $L^2$ norm of $q$ on $A$. So if $T-R$ is small enough we obtain $q=0$ on $A$; then one combines a unique continuation argument with a repeated application of the above to prove that $q=0$ on $A$ no matter what the $T$.

In the multidimensional case the above argument breaks down because of the angular Laplacian in (\ref{eq:vde}); all other parts of the argument work as in the one dimensional case. To carry out the above procedure we will need two estimates.
The first is a standard energy estimate for the wave equation and the second is an imitation of a sideways energy estimate for a one dimensional wave equation in $r,t$ where the roles of $r$ and $t$ are reversed.

\subsection{Energy identities}

\begin{figure}[h]
\begin{center}
\epsfig{file=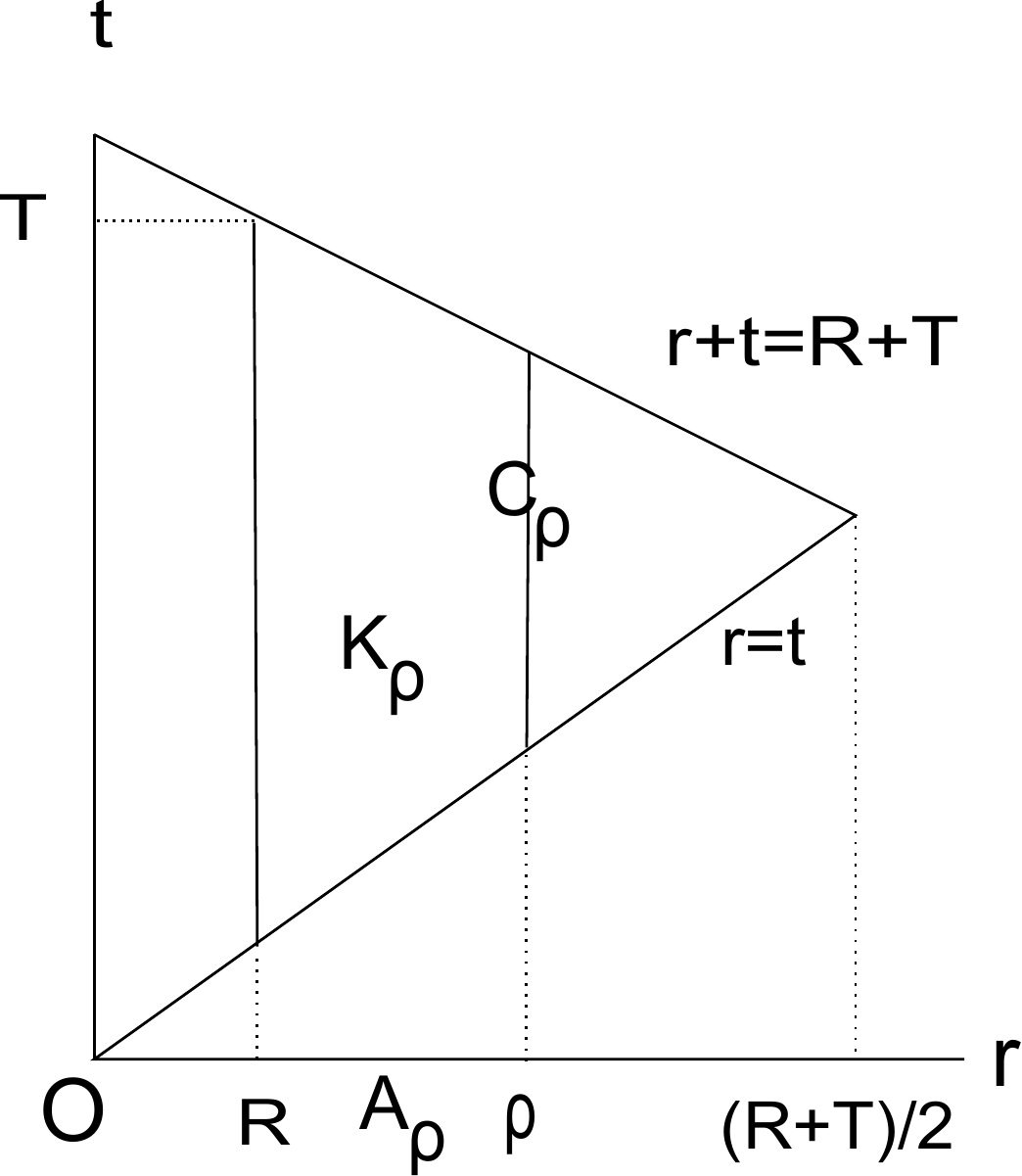, height=2.2in}
\end{center}
\caption{Sideways energy estimates}
\label{fig:sideways}
\end{figure}

For each $\rho \in [R, (R+T)/2]$, define (see Figure \ref{fig:sideways}) the sub-region 
\[
K_{\rho} :=  \{ (x,t) \in \R^3 \times \R \, : \, R \leq |x| \leq \rho, ~ |x| \leq t \leq R+ T-|x| \},
\]
the annular region
\[
A_\rho := \{ x \in \R^3 \, : \, R \leq |x| \leq \rho \},
\]
the vertical cylinder
\[
C_\rho :=  \{ (x,t) \in \R^3 \times \R \, : \, |x|=\rho, ~ \rho \leq t \leq R+T-r \},
\]
and for any function $w(x,t)$ let $\wb$ and $\wbb$ be the the restrictions of $w$ to the lower and upper characteristic cones, that is
\[
\wb(x) = w(x,|x|), \qquad  \wbb(x) = w(x,R+T-|x|).
\]

We derive some relations which lead to the estimates we need. These relations are either the standard energy identity or a sideways version of it. Suppose $w(x,t)$ satisfies
\[
w_{tt} - w_{rr} - \frac{1}{r^2} \Omega w  = F(x,t), \qquad (x,t) \in K.
\]
Define the ``sideways'' energy (we will assume a sum over $1\leq i < j \leq 3$)
\begin{align*}
 J(\rho) & := \int_{C_\rho} r^{-2}( w^2 + w_t^2 + |\nabla w|^2) \, dS_{x,t}
= \int_{C_\rho} r^{-2} (w_t^2 + w_r^2 + w^2 + r^{-2} (\Omega_{ij} w)^2 ) \, dS_{x,t}
\\
&= \int_{\rho}^{R+T-\rho} \int_S (w_t^2 + w_r^2 + w^2 + r^{-2} (\Omega_{ij}w)^2) (\rho \theta,t) \, d \theta \, dt.
\end{align*}
Multiplying the identity 
\begin{align}
2 w_r( w_{tt} -  w_{rr} - & r^{-2} \Omega w -  w)  - 4 r^{-2} \Omega_{ij} w_r \, \Omega_{ij} w +2r^{-3} (\Omega_{ij} w)^2
\nn
\\
&=  - (w_t^2 + w_r^2 + r^{-2} (\Omega_{ij} w)^2 +  w^2)_r + 2 (w_r w_t)_t - 2 \Omega_{ij}  (r^{-2} w_r \Omega_{ij} w )
\label{eq:wr}
\end{align}
by $r^{-2}$, integrating over the region $K_\rho$, using (\ref{eq:stokes}) and Stokes's theorem on a region in the $r,t$ plane, 
we obtain
\begin{align*}
\int_{K_\rho} r^{-2} &  \left ( 2 w_r(F-w)  - 4 r^{-2} \Omega_{ij} w_r \, \Omega_{ij} w +2r^{-3} (\Omega_{ij} w)^2 \right )
\; dx \, dt
\\
& = \int_S  \int_R^\rho \int_{r}^{R+T-r}  - (w_t^2 + w_r^2  +  r^{-2} (\Omega_{ij} w)^2 + w^2)_r + 
2 (w_r w_t)_t \, dt \, dr \, d \theta
\\
&= \int_S \int_R^T (w_t^2 + w_r^2 + r^{-2} (\Omega_{ij} w)^2 + w^2)(R \theta,t) \, dt \, d \theta
\\
& ~~~ - \int_S \int_\rho^{R+T-\rho} (w_t^2 + w_r^2 +  r^{-2} (\Omega_{ij} w)^2 + w^2)(\rho \theta,t) \, dt \, d \theta
\\
& ~~~
- \int_S  \int_R^\rho (w_t^2 + w_r^2 + r^{-2} (\Omega_{ij} w)^2 + w^2 - 2w_r w_t)(r \theta, R+T-r) \, dr \, d \theta
\\
& ~~~
-  \int_S  \int_R^\rho (w_t^2 + w_r^2 +  r^{-2} (\Omega_{ij} w)^2 + w^2 + 2w_r w_t)(r \theta, r) \, dr \, d \theta
\\
&= J(R) - J(\rho) - \int_{A_\rho} r^{-2} ( \wbb_r^2 + r^{-2} ( \Omega_{ij} \wbb)^2 + \wbb^2)(x) \, dx 
\\
& \qquad -
 \int_{A_\rho} r^{-2} ( \wb_r^2 +  r^{-2} ( \Omega_{ij} \wb)^2+ \wb^2)(x) \, dx.
\end{align*}
Hence
\begin{align}
J(\rho) + &  \int_{A_\rho} r^{-2} ( |\nabla \wbb|^2 +  \wbb^2)(x) \, dx +
 \int_{A_\rho} r^{-2} ( |\nabla \wb|^2 + \wb^2)(x) \, dx
+ \int_{K_\rho} 2r^{-5} (\Omega_{ij} w)^2 \, dx \, dt
\nn
\\
&= J(R)+ \int_{K_\rho} r^{-2} \left ( 2ww_r + 4r^{-2} \Omega_{ij} w_r \Omega_{ij}w - 2 F w_r \right )
\, dx \, dt,
\qquad R \leq \rho \leq \frac{R+T}{2}.
\label{eq:Jiden1}
\end{align}
This is the sideways energy identity we need.

\begin{figure}[h]
\begin{center}
\epsfig{file=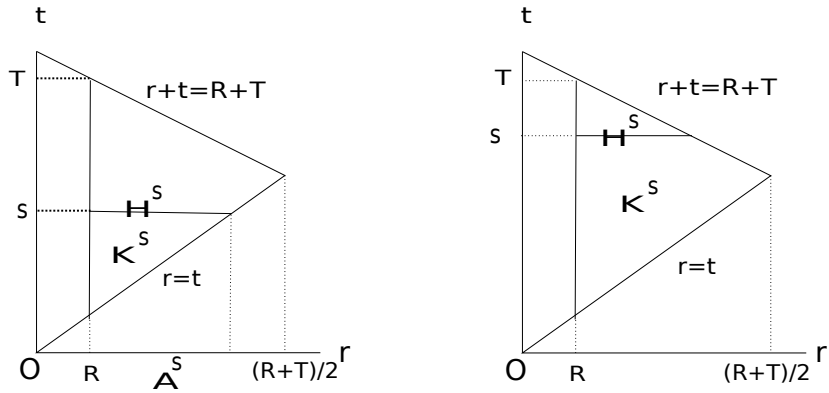, height=2in}
\end{center}
\caption{Standard energy estimate}
\label{fig:vertical}
\end{figure}
Next we derive the standard energy identity for the wave equation.
For any $s \in [R,T]$, define (see Figure \ref{fig:vertical}) the domain
\[
K^s = K \cap \{ (x,t) \in \R^3 \times \R \, : \, R \leq t \leq s \},
\]
$H^s$ the horizontal disk obtained by intersecting $K$ with the plane $t=s$, that is
\[
H^s = K \cap \{ t=s\},
\]
whose projection onto the plane $t=0$ is the annular region
\[
A^s := \{ x \in \R^3 \, : \, R \leq |x| \leq \min(s, R+T-s) \}.
\]
Next, we define the ``energy at time $s$'' for every $s \in [R,T]$ - the definition depends on $s \leq (R+T)/2$ or not because the geometry changes - see Figure \ref{fig:vertical}. For $ s \in [R, (R+T)/2]$, we define (summation over $1 \leq k< l \leq 3$)
\begin{align*}
E(s) &:= \int_{A^s} r^{-2} (w^2 + w_t^2+ w_r^2 + r^{-2} (\Omega_{kl} w)^2 ) (x,s) \, dx
\end{align*}
and for $s \in [(R+T)/2, T]$ we define
\begin{align*}
E(s) &:=  \int_{A^s} r^{-2} ( w^2 + w_t^2 + |\nabla w|^2)(x,s) \, dx + 
 \int_{R+T-s \leq |x| \leq (R+T)/2} r^{-2} ( \wbb(x)^2 + |(\nabla \wbb)(x)|^2 )\, dx.
\end{align*}

First take  $s \leq (R+T)/2$; multiplying the identity
\beqn
2w_t(w_{tt} - w_{rr} - r^{-2} \Omega w + w) = (w^2 + w_t^2 + w_r^2 + r^{-2} (\Omega_{kl} w)^2 )_t 
- 2(w_t w_r)_r - 2 r^{-2} \Omega_{kl} ( w_t \Omega_{kl} w)
\label{eq:wt}
\eeqn
by $r^{-2}$, integrating over the region $K^s$, and using (\ref{eq:divergence}), we obtain
\begin{align*}
\int_{K^s} & 2 r^{-2}  w_t ( F + w) \, dx \, dt
\\
& = E(s) +  2 \int_R^s \int_{|x|=R} r^{-2}  (w_t w_r)(x, s) \, dS_x \, dt
\\
& \qquad
- \int_{A^s} r^{-2} (w^2 + w_t^2 + w_r^2 + 2 w_t w_r + r^{-2} (\Omega_{kl} w)^2 )(x,|x|)
\, dx
\\
&=  E(s) +  2 \int_R^s \int_{|x|=R} r^{-2} (w_t w_r)(x, s) \, dS_x \, dt
- \int_{A^s} r^{-2} ( \wb(x)^2 + |(\nabla \wb)(x)|^2 ) \, dx.
\end{align*}
Next take $s \in [(R+T)/2, T]$; multiplying (\ref{eq:wt}) by $r^{-2}$, integrating over the region $K^s$,
using (\ref{eq:divergence}) we obtain
\begin{align*}
\int_{K^s} & 2 r^{-2} w_t  ( F + w) \, dx \, dt
\\
&= \int_{H^s} r^{-2} ( w^2 + w_t^2 + |\nabla w|^2) \, dx
\\
& ~~~~+ \int_{R+T-s \leq |x| \leq (R+T)/2} r^{-2} ( w^2 + w_t^2 + w_r^2 - 
2 w_t w_r + r^{-2} (\Omega_{kl} w)^2 )(x, R+T-|x|) \, dx
\\
& ~~~~-  \int_A r^{-2} (w_t^2 + w_r^2 + 2 w_t w_r + r^{-2} (\Omega_{kl} w)^2)(x,|x|) \, dx
+ 2\int_R^s \int_{|x|=R} r^{-2} (w_t w_r )(x,s) \, dS_x \, dt
\\
& =  \int_{H^s} r^{-2} ( w^2 + w_t^2 + |\nabla w|^2) \, dx 
+  \int_{R+T-s \leq |x| \leq (R+T)/2} r^{-2} ( \wbb(x)^2 + |(\nabla \wbb)(x)|^2 )\, dx
\\
& ~~~~ -  \int_A r^{-2}( \wb(x)^2 + |(\nabla \wb)(x)|^2) \, dx  + 
2 \int_R^s \int_{|x|=R} r^{-2} (w_t w_r )(x,t) \, dS_x \, dt
\\
&= E(s)  -  \int_A r^{-2}( \wb(x)^2 + |(\nabla \wb)(x)|^2) \, dx  + 
2 \int_R^s \int_{|x|=R} r^{-2} (w_t w_r )(x,t) \, dS_x \, dt.
\end{align*}
Hence, in either case, that is for any $s \in [R,T]$, we have
\begin{align}
E(s)  & \leq \int_A  r^{-2}( \wb(x)^2 + |(\nabla \wb)(x)|^2) \, dx  
+ 2 \int_{K^s} r^{-2} w_t  ( F + w)  \, dx \, dt
+  \int_C r^{-2} (w_t^2 + w_r^2) \, dS_x \, dt.
\label{eq:Eiden}
\end{align}

\subsection{Uniqueness}
We now show that if $v$ and $v_r$ are zero on $C$ then $q=0$ on $A$.
We apply (\ref{eq:Jiden1}) to $v=v_1 - v_2$; note that $ F = v_{tt} - v_{rr} - r^{-2} \Omega v = q_1v + qv_2$ and 
 $J(R)=0$ because the Cauchy data of $v$ is zero on $C$. Hence
\begin{align*}
J(\rho) + \int_{A_\rho} r^{-2} ( \vb^2 + |\nabla \vb|^2) &\leq
\int_{K_\rho} r^{-2} ( v v_r + 4 r^{-2} \Omega_{ij} v_r \Omega_{ij} v - 2v_r ( q_1 v + q v_2)
\\
&\cleq  \int_{K_\rho} r^{-2} (v^2 + v_r^2 + r^{-2}(\Omega_{ij} v)^2 + q^2 + r^{-2} (\Omega_{ij} v_r)^2 )
\\
&=  \int_R^\rho J(r) \, dr + \int_{A_\rho} r^{-2} q^2(x)  \left ( \int_r^{R+T-r} dt \right ) dx
+ \int_{K_\rho} r^{-4} (\Omega_{ij} v_r)^2
\\
& \leq
  \int_R^\rho J(r) \, dr + (T-R)  \int_A r^{-2} q^2(x) \, dx +   \int_K r^{-4} (\Omega_{ij} v_r)^2
\end{align*}
with the constant associated to $\cleq$ being $c_1 = 4\max(1, \|q_1\|_{L^\infty(A)}, \|v_2\|_{L^\infty(K)})$.
Hence, by Gronwall's inequality
\begin{align}
J(\rho) +  \int_{A_\rho} r^{-2} ( |\nabla p|^2 + p^2 )
\cleq (T-R)  \int_A r^{-2} q^2(x) \, dx +   \int_K r^{-4} (\Omega_{ij} v_r)^2,
\qquad R \leq \rho \leq \frac{R+T}{2},
\label{eq:JJtemp}
\end{align}
with the constant being $c_2 =c_1e^{c_1(T-R)}$. In particular
\beqn
J(\rho) \cleq 
  (T-R)  \int_A r^{-2} q^2(x) \, dx +   \int_K r^{-4} (\Omega_{ij} v_r)^2,
\qquad R \leq \rho \leq \frac{R+T}{2},
\label{eq:tempJ}
\eeqn
and taking $\rho = (R+T)/2$ in (\ref{eq:JJtemp}) we have
\beqn
 \int_A r^{-2} (  |\nabla p|^2 + p^2 ) \cleq 
 (T-R)  \int_A r^{-2} q^2(x) \, dx +   \int_K r^{-4} (\Omega_{ij} v_r)^2
\label{eq:J1}
\eeqn
with the constant $c_2$.
Integrating (\ref{eq:tempJ}) w.r.t $\rho$ over $[R, (R+T)/2]$ we obtain
\beqn
\int_K r^{-2} (v^2 + v_t^2+ |\nabla v|^2)
\cleq (T-R)^2 \int_A r^{-2} q^2(x) \, dx
+ (T-R)  \int_K r^{-4} (\Omega_{ij} v_r)^2.
\label{eq:J2}
\eeqn
So we can combine (\ref{eq:J1}), (\ref{eq:J2}) into
\begin{align}
\int_K  r^{-2} (v^2 + v_t^2 + |\nabla v|^2 ) &
+ 
 \int_A r^{-2} (  p^2 + |\nabla p|^2 )
 \cleq 
(T-R)  \int_A r^{-2} q^2 +   \int_K r^{-4} (\Omega_{ij} v_r)^2
\label{eq:Jest}
\end{align}
with the constant being $c_3 =(1+ T-R)c_2$.

The equation (\ref{eq:Jest}) would have been enough to prove Theorem 1 in the one dimensional case, because 
$|\nabla p|^2 \geq p_r^2 = q^2$ and the last term in (\ref{eq:Jest}) would not be there. Then by taking $T-R$ small enough we could have absorbed the second term on the RHS of (\ref{eq:Jest}) into the LHS and we would have proved the theorem for $T$ close to $R$. Then a unique continuation argument would prove the theorem for all $T>R$. However, in the three dimensional case we do have the last term in (\ref{eq:Jest}) which cannot be absorbed in the LHS because it involves second order derivatives of $v$ - we will estimate it in terms of $p$ using the standard energy estimate for the wave operator.

Fix an $i,j$ pair with $i<j$. We apply (\ref{eq:Eiden}) to the function $w = \Omega_{ij} v$, noting that
$\Omega_{ij}$ commutes with $\Omega$. Note that from (\ref{eq:vde}) and (\ref{eq:vic}) we have
\[
w_{tt} - w_{rr} - \frac{1}{r^2} \Omega w  = F
\]
with
\beqn
F(x,t) := q_1 w +  (\Omega_{ij} q_1)v +  (\Omega_{ij} q)v_2 + q \Omega_{ij} v_2.
\label{eq:Mdef}
\eeqn
and
\beqn
\wb(x,|x|) = (\Omega_{ij} p)(x).
\eeqn
Further, since the Cauchy data of $v$ is zero on $C$, so the Cauchy data of $w$ is zero on $C$.
Hence from (\ref{eq:Eiden}) we have 
\begin{align*}
E(s) & \leq \int_A r^{-2} ( (\Omega_{ij}p)^2 + |\nabla \Omega_{ij} p|^2)
+ \int_{K^s} r^{-2} ( w^2 + w_t^2 + F^2)
\\
& \cleq  \int_A r^{-2}  ( (\Omega_{ij}p)^2 + |\nabla \Omega_{ij} p|^2)
+ \int_{K^s} r^{-2} (w^2 + w_t^2 + v^2 + q^2 + (\Omega_{ij} q)^2 )
\\
& \cleq \int_R^s E(t) \, dt + 
\int_A  r^{-2} ( p^2 + |\nabla p|^2  + |\nabla \Omega_{ij}p|^2 )
+ \int_K r^{-2} v^2 
\end{align*}
with the constant being $c_4 = 2 \max( 1, (R+T)^2, \|q_1\|_\infty, \|\Omega_{ij} q_1 \|_\infty, \|v_2\|_\infty )$.
So from Gronwall's inequality we have 
\begin{align*}
E(s) \cleq \int_A  r^{-2} ( p^2 + |\nabla p|^2 +  |\nabla \Omega_{ij} p|^2  )
+ \int_K r^{-2} v^2 ,
\qquad R \leq s \leq T
\end{align*}
with the constant being $c_5 = c_4 e^{c_4 (T-R)}$.
Integrating this w.r.t $s$ over the interval $[R,T]$ we obtain
\begin{align*}
\int_K r^{-2}  ( w^2 + w_t^2 + |\nabla w|^2)
 \leq c_5(T-R) \left (
\int_A  r^{-2} ( p^2 +  |\nabla p|^2 + |\nabla \Omega_{ij}p|^2  )
+ \int_K r^{-2} v^2
\right );
\end{align*}
hence, since $w = \Omega_{ij} v$,
\begin{align}
\int_K r^{-4} & (\Omega_{ij} v_r)^2
\leq c_5R^{-2}(T-R) \left (
\int_A  r^{-2} ( p^2 + |\nabla p|^2 + |\nabla  \Omega_{ij}p|^2  )
+ \int_K r^{-2} v^2
\right ).
\label{eq:Eest}
\end{align}
Using this in (\ref{eq:Jest}), we have
\begin{align}
\int_K  r^{-2} v^2  &
+ 
 \int_A r^{-2} (p^2 + |\nabla p|^2) 
\nn
\\
&
\cleq 
(T-R)  \int_A r^{-2} ( p^2 + |\nabla p|^2 + | \nabla \Omega_{ij} p|^2 )
+ (T-R) \int_K r^{-2} v^2
\label{eq:tempvw}
\end{align}
with the constant $c_6 = \max(c_3, c_3c_5 R^{-2})$. However, $q$ is in $Q_\gamma$ so
\begin{align*}
\int_A r^{-2} |\nabla (\Omega_{ij} p)(x)|^2 \, dx
&=
\int_R^{(R+T)/2} \int_{| \theta|=1} (\nabla \Omega_{ij} p)(r \theta) ^2 \, d \theta \, dr
\\
& \leq
\gamma  \int_R^{(R+T)/2} r^2 \int_{| \theta|=1} (p^2 + | \nabla p|^2)(r \theta) \, d \theta \, dr
\\
& \leq \gamma (R+T)^2 \int_A r^{-2} (p^2 + |\nabla p|^2).
\end{align*}
Using this in (\ref{eq:tempvw}),  we see that $p=0$ on $A$ if $T-R$ is small enough - depending on $\gamma$, $c_6$ and $R+T$. Now $v(x,|x|) = p(x)$ and $v=0$ on $|x|=R$ so $p=0$ on $|x|=R$,
that is $ \int_0^R q( \sigma \theta) \, d \sigma=0$ for all unit vectors $\theta$. Hence
\[
\int_R^r q( \sigma \theta) \, d \sigma =0, \qquad R \leq r \leq T
\]
which implies $q(x)=0$ when $R \leq |x| \leq T$, provided $T-R$ is small enough.

Actually, adjusting the height of the downward pointing cone, what we have shown is the 
following: there is a $\delta>0$ dependent only on
$\gamma, R, T, \|q_1\|_{C^1(A)}, \|v_2\|_{C^1(K)}$, so that if, for some $R^* \in [R, (R+T)/2]$,
$v$ and $v_r$ are zero on the cylinder 
\[
\{ (x,t) \; : \; |x|=R^*, ~ R^* \leq t \leq R^* + 2 \delta \},
\]
then $q=0$ on $R^* \leq |x| \leq R^* + \delta$, with the obvious modification in the assertion if $R^*+ \delta > (R+T)/2$.
We use this observation to prove that $q=0$ for any $R,T$.

Since $v$ and $v_r$ are zero on $C$, then from the above claim, we have $q=0$ on $R \leq |x| \leq R+ \delta$. Let $u=u_1-u_2$ where $u_1, u_2$ are solutions to (\ref{eq:ude}), (\ref{eq:uic}) for $q=q_1, q_2$. Then, $u$ satisfies the homogeneous equation
\[
u_{tt} - \Delta u - q_1 u = 0
\]
over the region $K_\rho$ where $\rho=R+\delta$. Now $u$ and $u_r$ are zero on $C$, and $q_1$ is independent of $t$, so by the Robbiano-Tataru unique continuation theorem (see Theorem 3.16 in [KKL01]) we have $u=0$ in the region $K_\rho$; in particular $u$ and $u_r$ are zero on $C_\rho$ and hence $v,v_r$ are zero on $C_\rho$. Now repeat the above argument, except $R$ is replaced by $R+\delta$; this argument repeated will complete the proof of Theorem \ref{thm:reflection}.


\section{Proof of Theorem \ref{thm:linearized}}

Let $\{ \phi_n(x) \}_{n=1}^\infty$ be a sequence of homogeneous harmonic polynomials on $\R^3$ so that their restrictions to the unit sphere $S$ form an orthonormal basis on $L^2(S)$ - see Chapter 4 of 
\cite{SW71}. Let $k(n)$ be the degree of homogeneity of $\phi_n$. Then $q(x)$ and $u(x,t)$ have spherical harmonic decompositions in $L^2(S)$ given by
\begin{align*}
q(r \theta) = \sum_{n=1}^\infty q_n(r) r^{k(n)} \phi_n(\theta),
\qquad
u(r \theta, t) =  \sum_{n=1}^\infty u_n(r, t) r^{k(n)} \phi_n(\theta)
\end{align*}
where
\[
r^{k(n)} q_n(r) = \int_{|\theta|=1} q(r \theta) \, \phi_n(\theta) \, d \theta,
\qquad
r^{k(n)} u_n(r,t) = \int_{|\theta|=1} u(r \theta,t) \, \phi_n(\theta) \, d \theta.
\]
Since $u$ and $q$ are smooth, we may 
show\footnote{Use the definition of $q_n$ and $u_n$, observe that the $\phi_n(\theta)$ are eigenvalues of the spherical Laplacian, and use the Divergence Theorem on $S$ to transfer the Laplacian from the $\phi_n$ to $q$ or $u$ - see Theorems 2 and 4 in \cite{See66}.} 
that $q_n(r)$ and $u_n(r,t)$ decay as $n^{-p}$ for large $n$ for any positive integer $p$, uniformly in $r,t$. Hence the series also converge in the $C^2$ norm.

To prove the theorem, it will be enough to prove that $q_n(r)=0$ on $R \leq r \leq (R+T)/2$ for all $n \geq 1$.
One may show that for sufficiently regular $f$ (see page 1235 of \cite{FPR04})
\[
\Delta \left ( f(r,t)  r^{k(n)} \phi_n(\theta) \right ) = r^{k(n)} \phi_n(\theta) (f_{tt} - f_{rr} - \frac{ 2 k(n) -2}{r} f_r )
\]
hence, using (\ref{eq:ulinde}), (\ref{eq:ulinic}), the $u_n(r,t)$ are solutions of the one dimensional Goursat problems
\begin{gather*}
\pa_t^2 u_{n} - \pa_r^2 u_{n,rr} - \frac{ 2 k(n) - 2}{r} \pa_r u_n - q_b u_n = q_n u_b, \qquad t \geq |r|
\label{eq:unde}
\\
u_n(r,|r|) = \int_0^1 \sigma^{ k(n)} q_n( \sigma r) \, d \sigma.
\label{eq:unic}
\end{gather*}
The hypothesis of the theorem implies that $u_n(R,t)$ and $(\pa_r u_n)(R,t)$ are zero for $R \leq t \leq T$. So repeating the standard argument for one dimensional hyperbolic inverse problems with reflection data,  as in \cite{Sym86}, or repeating just the sideways energy argument in the proof of Theorem \ref{thm:reflection} without the complication of the angular terms, one may show that $q_n(r)=0$ for $R \leq r \leq (R+T)/2$.

\section{Proof of Theorem \ref{thm:transmission}}\label{sec:transmission}

\begin{figure}[h]
\begin{center}
\epsfig{file=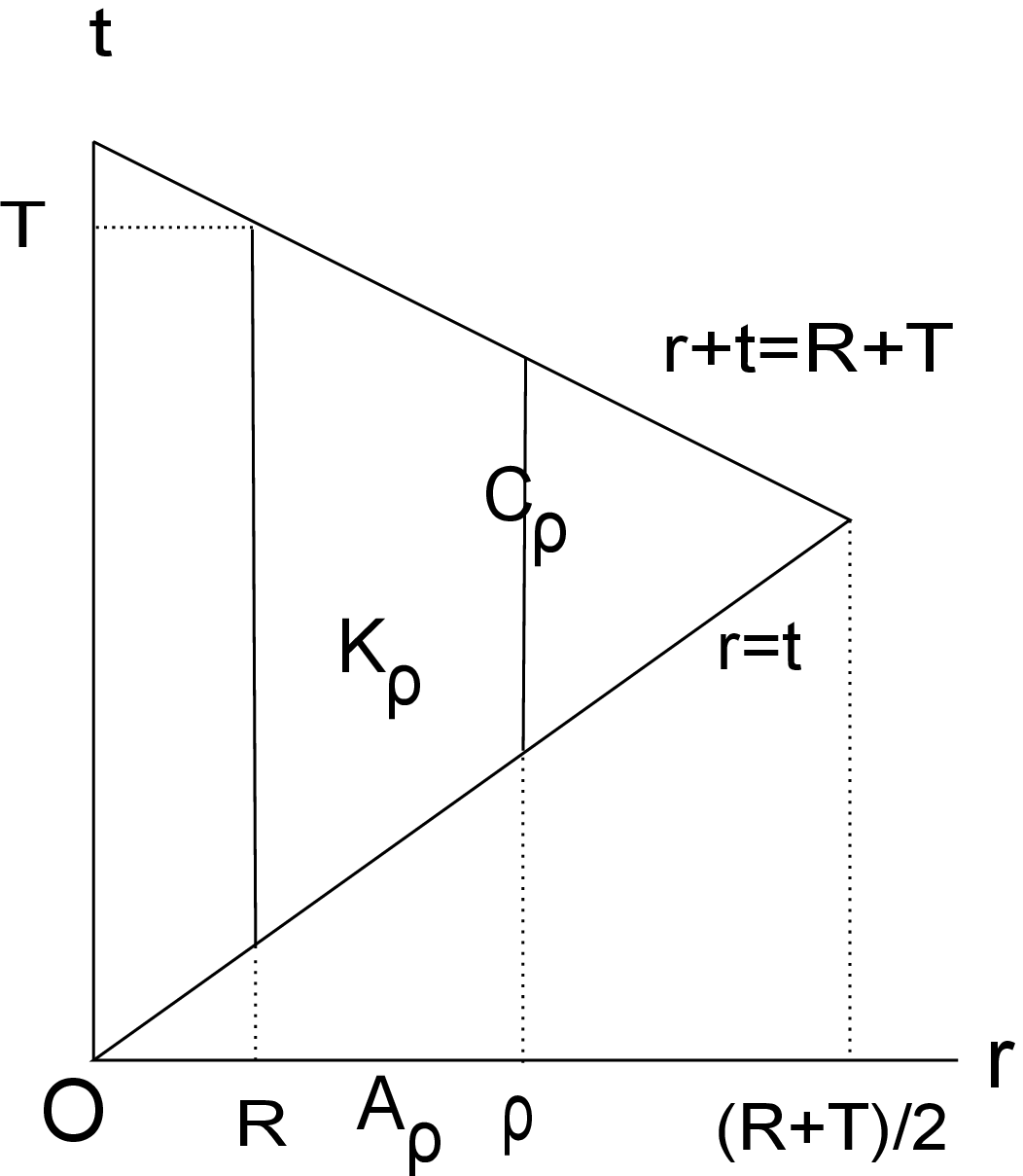, height=2in}
\end{center}
\caption{Transmission data problem}
\label{fig:small}
\end{figure}
Let (see Figure \ref{fig:small}) $B$ denote the origin centered ball of radius $R$ in $\R^3$,  $D$ the region
\[
D := \{ (x,t) \in \R^3 \times \R \, : |x| \leq R, ~ |x| \leq t \leq T \},
\]
and as before $C$ the cylinder
\[
C := \{ (x,t) \,: \, |x|=R, ~ R \leq t \leq T \}.
\]
Let $u_i$, $i=1,2$ be the solutions of (\ref{eq:ude}), (\ref{eq:uic}) when $q=q_i$; define $q=q_1-q_2$ and $u=u_1-u_2$. Then $u$ satisfies
\begin{gather}
u_{tt} - \Delta u - q_1u = qu_2, \qquad (x,t) \in D
\label{eq:udetran}
\\
u(x,|x|) = \int_0^1 q(\sigma x) \, d \sigma.
\label{eq:uictran}
\end{gather}
Then, restricting attention to the cylindrical region $B \times [R,T]$, from \cite{KM91} we have the following stability estimate for the time-like Cauchy problem (note $T>3R$): there is a constant $C_1$ dependent only on $M,R,T$ so that
\beqn
\| u(\cdot,t) \|_{H^1(B)}^2 +  \| u_t(\cdot,t) \|_{L^2(B)}^2
\leq C_1 \left ( \|q u_2\|_{L^2(B \times [R,T] )}^2 +  \|u\|_{H^1(C)}^2 + \|u_r \|_{L^2(C)}^2 \right ),
\qquad R \leq t \leq T.
\label{eq:tem1}
\eeqn
Next, if we multiply (\ref{eq:udetran}) by $u_t$ and use the techniques for standard energy estimates (backward in time) on the region $|x| \leq t \leq R$, we obtain
\beqn
\int_B |\ub(x)|^2 + |\nabla \ub(x)|^2 \, dx \leq C_2 \left ( \iint_{|x| \leq t \leq R} |q u_2|^2 \, dx \, dt +
\| u(\cdot,R) \|_{H^1(B)}^2 +  \| u_t(\cdot,R) \|_{L^2(B)}^2 \right )
\label{eq:tem2}
\eeqn
where $\ub(x) = u(x,|x|)$ and $C_2$ depends only on $M,R$.
Hence, combining (\ref{eq:tem1}), (\ref{eq:tem2}) we obtain
\beqn
\int_B |\ub(x)|^2 + |\nabla \ub(x)|^2 \, dx
\leq C_3 \left ( \|q u_2 \|_{L^2(D)}^2 + 
 \|u\|_{H^1(C)}^2 + \|u_r \|_{L^2(C)}^2 \right )
\eeqn
where $C_3$ depends only on $R,T,M$.
Now $ r \ub(x) = \int_0^r q(s \theta) \, ds$, hence
$ q(x) = (r \ub)_r = \ub + r \ub_r$. So
\[
q^2 \leq 2 (\ub^2 + r^2 \ub_r^2) \leq 2 \max(1, R^2) ( \ub^2 + \ub_r^2)
\leq 2 \max (1, R^2) ( \ub^2 + |\nabla \ub |^2 ),
\]
and
\beqn
\|q\|_{L^2(B)}^2 \leq C_4 \left ( \|q u_2 \|_{L^2(D)}^2 + 
 \|u\|_{H^1(C)}^2 + \|u_r \|_{L^2(C)}^2 \right )
\label{eq:final1}
\eeqn
with $C_4$ dependent only on $R,T,M$. Finally, using Theorem \ref{thm:existence}, we have
\[
\| q u_2 \|_{L^2(D)} \leq \|u_2\|_{L^\infty(D)} \, \|q\|_{L^2(D)}
\leq \N(T,\|q_2\|_\infty) \, \|q\|_{L^2(D)}
\]
where the $\|q_2\|_\infty$ norm is over the region $|x| \leq (R+T)/2$. Since $ \N(T,\|q_2\|_\infty)$ goes to zero
as $\|q_2\|_\infty$ approaches $0$, we can choose a $\delta>0$ so that 
\[
C_4  \N(T,\|q_2\|_\infty) < \frac{1}{2}
\]
if $ \| q_2 \|_\infty \leq \delta$; note that this $\delta$ will depend only on $R,T,M$. Using this in (\ref{eq:final1}), we conclude that if $\|q_2\|_\infty \leq \delta $ then
\beqn
\|q\|_{L^2(B)}^2 \leq C_5 \left (
 \|u\|_{H^1(C)}^2 + \|u_r \|_{L^2(C)}^2 \right )
\label{eq:final}
\eeqn
with $C_5$ dependent only on $R,T,M$.

\bibliography{REFS}{}
\bibliographystyle{abbrv}

\end{document}